\documentclass[10 pt, journal]{IEEEtran}
\IEEEoverridecommandlockouts                               % This command is only
\usepackage{graphicx} % for pdf, bitmapped graphics files
\usepackage{epsfig} % for postscript graphics files
\usepackage{mathptmx} % assumes new font selection scheme installed
\usepackage{times} % assumes new font selection scheme installed
\usepackage{amsmath} % assumes amsmath package installed
\usepackage{amssymb}  % assumes amsmath package installed

\title{Design and Implementation of Repetitive Control based Noncausal Zero-Phase Iterative Learning Control}

\author{K.~Krishnamoorthy
        and~Tsu-Chin~Tsao,~\IEEEmembership{Senior Member,~IEEE}% <-this % stops a space
\thanks{This work was supported in part by the National Science Foundation under grants DMI 0327077 and CMMI 0751621. This paper was presented in parts at the ASME Mechanical Engineering Congress and Exposition, Anaheim, CA, Nov 2004 and IEEE Conference on Decision and Control, Atlantis, Bahamas, Dec 2004}
\thanks{K.~Krishnamoorthy was with the University of California at Los Angeles, CA 90095. He is now with the Air Force Research Laboratory/ Control Design and Analysis Branch, Wright-Patterson Air Force Base, OH 45433. Email: krishnak@ucla.edu}
\thanks{T.-C.~Tsao is with the Mechanical and Aerospace Engineering Department, University of California at Los Angeles, CA 90095. Email: ttsao@ucla.edu}
}

\begin{document}

\maketitle

\thispagestyle{empty}
\pagestyle{empty}

%%%%%%%%%%%%%%%%%%%%%%%%%%%%%%%%%%%%%%%%%%%%%%%%%%%%%%%%%%%%%%%%%%%%%%%%%%%%%%%%
\begin{abstract}
Discrete-time domain Iterative Learning Control (ILC) schemes inspired by Repetitive control algorithms are proposed and analyzed. The well known relation between a discrete-time plant (filter) and its Markov Toeplitz matrix representation has been exploited in previous ILC
literature. However, this connection breaks down when the filters
have noncausal components. In this paper we provide a formal
representation and analysis that recover the connections between
noncausal filters and Toeplitz matrices. This tool is then applied to
translate the anti-causal zero-phase-error prototype repetitive control scheme to the design of a stable and fast converging ILC algorithm. The learning gain matrices are
chosen such that the resulting state transition matrix has a \textit{Symmetric Banded Toeplitz} (SBT) structure. It is shown that the well known sufficient condition for repetitive control closed loop stability based on a filter's frequency domain $H_\infty$ norm is also sufficient for ILC convergence and that the condition becomes necessary as the data length approaches infinity. Thus the ILC learning
matrix design can be translated to repetitive control loop shaping
filter design in the frequency domain.
\end{abstract}

%%%%%%%%%%%%%%%%%%%%%%%%%%%%%%%%%%%%%%%%%%%%%%%%%%%%%%%%%%%%%%%%%%%%%%%%%%%%%%%%
\section*{Nomenclature}
\label{sec:nom}
\noindent
In this paper upper case symbols (e.g., $G$) represent discrete transfer functions.
Upper case bold symbols (e.g., $\mathbf G$) represent corresponding matrices and
lower case bold symbols (e.g, $\mathbf u$) represent vectors in the ``lifted'' domain. Lower case symbols (e.g, $b$) represent scalar values. Also $\mathbf 0$ and $\mathbf I$ stand for the zero and identity matrices respectively. In general $\mathbf M_{a,b}$ stands for a matrix of dimensions $a\times b$.
%We assume all the signals encountered hereafter belong to $L_2[0,T]$ where $T$ is finite.

%%%%%%%%%%%%%%%%%%%%%%%%%%%%%%%%%%%%%%%%%%%%%%%%%%%%%%%%%%%%%%%%%%%%%%%%%%%%%%%%%%%
\section{Introduction}
\label{intro}
\noindent
Iterative Learning Control is a common methodology used in reducing
tracking errors trial-by-trial for systems that operate repetitively
\cite{Ari,Kevb,Bien}.
In such systems, the reference is usually unchanged from one iteration to the next.
An expository view of common ILC schemes is provided in \cite{Kev98,Kev07}. General
convergence conditions based on operator norm for different ILC schemes can be found
in \cite{Kev89}. The relationship between learning filter design in the ``lifted'' domain and causal filter design in the frequency domain is detailed in \cite{Owe92,Owe96,Dim,Nor02}.
We wish to extend this relationship to noncausal filter design in the time axis. Higher order ILC schemes both in time and iteration axis are described in
\cite{Kev,Kev02,Kev03}.The learning filter design in \cite{Kev03} is based on minimization of the tracking error norm. For monotonic convergence, high-order ILC schemes in iteration axis alone may not suffice \cite{Nor,Jian}.

Causal ILC laws result in the transition matrix having a lower triangular Toeplitz structure. The Toeplitz matrix structure and it's commutability is lost when we deal with noncausal filters. In this paper, a treatment for this difficulty is proposed and also used to facilitate the ILC design. This design is analogous to noncausal zero-phase error compensation and filtering used in repetitive control design\cite{Tom85}. Connections between repetitive control and ILC are well established in the literature \cite{Long00,Long02, Long03, Xu06, Okko, Chen}. Zero-phase based ILC laws are known to result in good transient behavior and have nice robustness properties. But thus far they have been developed for infinite time signals that can only be applied to sufficiently long finite time signals \cite{Long94a}. Using noncausal operators in ILC design is in itself not a new idea \cite{Owe98,Verw}. As reported in \cite{Verw,Verw05}, there is good reason to consider noncausal operators especially for non-minimum phase systems and also for minimum phase systems with high relative degree. In particular, zero-phase noncausal filters have been employed as learning gains for servo-positioning applications \cite{Wu09,Leang,Deva05}. But most of the earlier ILC work on noncausal filters make use of infinite-time domain (continuous or discrete) theory and employ related frequency domain conditions to arrive at conditions for iteration convergence. It is usually (tacitly) assumed that the same conditions will hold when the learning controller is implemented using discrete finite-time signals. In this context, we point out some crucial observations made in \cite{Long00}:
\begin{enumerate}
\item[1)] Stability condition for ILC based on steady state frequency response strictly applies only to infinite-time signals. 
\item[2)] Conclusions based on frequency domain conditions apply only to parts of the trajectory for which steady state frequency response describes the input-output relation.
\end{enumerate}
For causal filters, the exact relationship between ILC convergence (true stability condition) and the frequency domain (approximate stability condition) are brought out in the same paper. 
%But as noted in \cite{Long00}, the frequency domain condition for ILC based on steady state frequency response strictly applies only to infinite-time signals. Furthermore, the relationship . One important observation in the paper \cite{Long00} is that conclusions based on frequency domain conditions apply only to parts of the trajectory for which steady state frequency response describes the input-output relation.
%One can in fact easily construct scenarios where the approximate stability condition is satisfied but the true stability condition is not (see example in Section~\ref{sec:ex}).
%as the data length $n\rightarrow\infty$, if the transition matrix is toeplitz, then the true stability condition for ILC convergence coincides with the approximate condition based on a steady state frequency response.
Similar arguments can be made for noncausal ILC as well and hence one cannot assume, in general, that the approximate stability condition is a sufficient condition for iteration error convergence. We therefore propose a special choice of the learning gain matrices that translates the ILC matrix design to repetitive control loop shaping design in the frequency domain. 
%In this paper, 
In particular, the noncausal learning gain matrices we choose result in a symmetric banded toeplitz (SBT) cycle-to-cycle transition matrix. This special structure enables us to use the frequency domain approximate stability condition as a sufficient condition for iteration convergence. Furthermore, we show that the approximate stability condition becomes both necessary and sufficient for iteration convergence when the data length approaches infinity. 
%This paper is a compilation of the results presented in \cite{Kri1,Kri2}. Also 
The scheme presented herein was successfully implemented on a dual stage fast tool servo system \cite{Kri3}.
%ILC laws motivated by stable plant inversion based feedforward control \cite{Tom85}, have not been considered.
\\\indent
The rest of this paper is organized as follows: Section~\ref{sec:ilc} introduces
common ILC based on causal compensations such as P or PD type ILC. Section~\ref{sec:rep} reviews the so called prototype repetitive control and motivates an analogous ILC scheme.
Section~\ref{sec:modrep} reviews the repetitive control algorithm modified for robustness and
the central theme, the modified repetitive ILC algorithm.
Section~\ref{sec:stab} provides a convergence analysis of the
proposed algorithm with sufficient conditions established in the frequency domain.
Section~\ref{sec:mon} establishes conditions for monotonic convergence of the error for a special choice of learning gain matrix. Section~\ref{sec:ex} gives a simple simulation example to bring out a key requirement for the said convergence. Section~\ref{sec:exp} details the experimental verification of the proposed scheme and a simulation which shows how the ILC scheme can improve over stable inversion based feedforward control.

%%%%%%%%%%%%%%%%%%%%%%%%%%%%%%%%%%%%%%%%%%%%%%%%%%%%%%%%%%%%%%%%%%%%%%%%%%%%%%%%%%%

\section{Problem Setup}
\label{sec:ilc}
\noindent
Let us define the learning control, plant output and desired output by supervectors $\mathbf u_k,\mathbf y_k$ and $\mathbf r$ respectively.
\begin{eqnarray}
\mathbf u_k &=& [\begin{array}{ccc}u(0)& \ldots &
                           u(n-1)\end{array}]^T \nonumber \\
\mathbf r &=& [\begin{array}{ccc}r(d)& \ldots &
                           r(n+d-1)\end{array}]^T \nonumber \\
\mathbf y_k &=& [\begin{array}{ccc}y(d)& \ldots &
                           y(n+d-1)\end{array}]^T
\end{eqnarray}
where $n$ is the length of each trial, $k$ is the iteration index
and $d$ is the relative degree of the plant.
%Without loss of generality we assume that $d=1$ from here on.
The plant output can then be represented in the so called ``lifted'' domain (introduced in \cite{Long88} as a mathematical framework to represent ILC systems)  as
\begin{eqnarray}
\label{eq:plant}
\mathbf y_k &=& \mathbf G\mathbf u_k
\end{eqnarray}
Where $\mathbf G$ is a lower triangular matrix of Markov parameters of the Linear Time Invariant (LTI) plant given by
\begin{eqnarray}
\label{eq:mark}
\left[\begin{array}{ccccc}
      h_d &0&0&\ldots&0\\
      h_{d+1} & h_{d} &0&\ldots&0\\
      h_{d+2} & h_{d+1} & h_{d} & \ldots & 0\\
      \vdots&\vdots&\vdots&\ddots &\vdots\\
      h_{n+d-1}& h_{n+d-2} & h_{n+d-3}& \ldots & h_d\\
\end{array}\right]
\end{eqnarray}
We assume the initial condition is the same for each iteration and hence it
does not appear in the error propagation. We shall
ignore it hereafter as is commonly done in the literature.
Now a generic ILC update law would be of the form
\begin{eqnarray}
\label{eq:ilc}
\mathbf u_{k+1} = \mathbf T_u\mathbf u_k + \mathbf T_e(\mathbf r-\mathbf y_k)
\end{eqnarray}
where $\mathbf T_u$ and $\mathbf T_e$ are square matrices. This leads to the
propagation equation
\begin{eqnarray}
\mathbf u_{k+1} &=& (\mathbf T_u-\mathbf T_e\mathbf G)\mathbf u_k + \mathbf T_e\mathbf r \nonumber \\\\
\mathbf e_{k+1} &=& \mathbf r-\mathbf G(\mathbf T_u\mathbf u_k+\mathbf T_e\mathbf e_k) \nonumber\\
        &=& (\mathbf T_u-\mathbf G\mathbf T_e)\mathbf e_k + (\mathbf I-\mathbf T_u)\mathbf r \nonumber \\
        && if \quad \mathbf G,\mathbf T_u \quad commute \nonumber
\end{eqnarray}
We note that two matrices commute if they are lower triangular Toeplitz or circulant \cite{Gray}
%%%%%%%%%%%%%%%%%%%%%%%%%%%%%%%%%%%%%%%%%%%%%%%%%%%%%%%%%%%%%%%%%%%%%%%%%%%%%%%%%%%
\subsection{Arimoto P-type ILC Algorithm}
\label{sec:ppd}
\noindent
First we define the tracking error as the vector
\begin{eqnarray}
\label{eq:errdef}
\mathbf e_k = \mathbf r-\mathbf y_k
\end{eqnarray}
The \textit{Arimoto P-type} learning law \cite{Ari} is given by
\begin{eqnarray}
u_{k+1}(t) = u_k(t) + \alpha e_k(t)
\end{eqnarray}
which we get by substituting $\mathbf T_u=\mathbf I$ and $\mathbf T_e=\alpha \mathbf I$
in (\ref{eq:ilc}).
The tracking error propagates according to
\begin{eqnarray}
\mathbf e_{k+1} &=& \mathbf r - \mathbf y_{k+1} \nonumber\\
&=& \mathbf r - \mathbf G\mathbf u_{k+1} \nonumber\\
&=& \mathbf r - \mathbf Gu_k - \alpha \mathbf G\mathbf e_k \nonumber \\
&=& (\mathbf I-\alpha \mathbf G)\mathbf e_k
\end{eqnarray}
A necessary and sufficient condition for iteration convergence is
that the state transition matrix $(\mathbf I-\alpha \mathbf G)$ be stable. Since the state transition matrix is
lower triangular this translates to $|1-\alpha h_d| < 1$.
If it is stable, we have the fixed point of the iteration given by
\begin{eqnarray}
\label{eq:arimo}
\mathbf e_\infty &=& 0 \nonumber \\
\mathbf u_\infty &=& \mathbf u_\infty + \alpha (\mathbf r-\mathbf G\mathbf u_\infty)
\nonumber \\
&=& \mathbf G^{-1}\mathbf r
\end{eqnarray}
The convergence condition can be easily met by an appropriate
choice of the learning gain $\alpha$. Although this is a necessary and sufficient
for iteration convergence, the convergence is not necessarily monotonic. To guarantee
monotonicity \cite{Kev01} we require
\begin{eqnarray}
\label{eq:mon1}
|1-\alpha h_d|+|\alpha|\sum_1^{n-1}|h_{d+i}| < 1
\end{eqnarray}
If the Markov parameter $h_d$ is close to zero, the matrix $\mathbf G^{-1}$  becomes
numerically unstable
(very large
entries) leading to undesirable control input $\mathbf u_\infty$. This could lead to actuator
saturation and hence may not be implementable in practise. We could instead have a
\textit{PD-type} learning law \cite{Kev03}
\[ u_{k+1}(t) = u_k(t) + \alpha e_k(t) + \beta e_k(t-1) \]
which would result in a \textit{2-band} $\mathbf T_e$ matrix with zeros above the main
diagonal.
We can extend the update law to being noncausal in time which
translates to $\mathbf T_e$ and possibly $\mathbf T_u$ matrices having non-zero elements
 above
the main diagonal. For example, the so called zero-phase filter of order $nq$
\begin{eqnarray}
\label{eq:qfil}
Q(z,z^{-1}) &=& q_0 + q_1(z+z^{-1}) + \ldots + q_{nq}(z^{nq}+z^{-nq}) \nonumber \\
\sum_{i=0}^{nq} q_i &=& 1
\end{eqnarray}
will result in a symmetric banded Toeplitz matrix. There exists a sufficient condition for iteration convergence in the frequency domain when using causal filters \cite{Owe96}.
We would like to extend this result to the noncausal update law as well.
It turns out that this transition is not straightforward and the reasons are
elucidated later.
If we were to restrict the noncausal filters to be of the form (\ref{eq:qfil})
the resulting state transition matrix has a favorable structure which enables us to derive
sufficient conditions for iteration convergence.

%%%%%%%%%%%%%%%%%%%%%%%%%%%%%%%%%%%%%%%%%%%%%%%%%%%%%%%%%%%%%%%%%%%%%%%%%%%%%%%%%%
\section{Prototype Repetitive Control}
\label{sec:rep}
\noindent
Repetitive control \cite{Tom89,Tom90}
is a special type of controller that handles periodic signals
based on the internal model principle \cite{Won75}. The repetitive control scheme rejects
disturbances appearing at a known fundamental frequency and its
harmonics.
Given the stable causal plant transfer function
\begin{eqnarray}
\label{eq:filt}
G(z^{-1}) &=& z^{-d}G^+(z^{-1})G^-(z^{-1})
\end{eqnarray}
where $d$ is the relative degree of the plant and we have split the
transfer function into invertible and non-invertible parts.
%The monic stable polynomials $B^+$ and $D$ are defined by
%\[
%\begin{array}{ccl}
%B^+(z^{-1}) & = & 1+b_1z^{-1}+b_2z^{-2}+\ldots+b_mz^{-m} \\
%D(z^{-1}) & = & 1+d_1z^{-1}+d_2z^{-2}+\ldots+d_pz^{-p}
%\end{array}
%\]
The prototype repetitive controller \cite{Tom90} based on stable plant
inversion is given by
\begin{eqnarray}
K_{rep}(z^{-1}) &=& \alpha\frac{z^{-N+d}G^-(z)}{b G^+(z^{-1})(1-z^{-N})} \nonumber \\
b &=& max_\omega|G^-(e^{-j\omega})|^2,\quad \omega\in[0,\pi]\
\end{eqnarray}
This can be written is feedback form as
\begin{eqnarray}
\label{eq:prepro}
u(t) = u(t-N) + \alpha\frac{\left[z^{-nu}G^-(z)\right]}{bG^+(z^{-1})} e(t-N+d+nu)
\end{eqnarray}
where $d$ is the relative degree of the plant,
$nu$ denotes the number of unstable zeros and $N$ is the period.
Introducing the variable $u'(t) = bG^+(z^{-1})u(t)$ we get
\begin{eqnarray}
\label{eq:proto}
u'(t) = u'(t-N) +  \alpha \left[z^{-nu}G^-(z)\right]e(t-N+d+nu)
\end{eqnarray}
Note the abuse of notation in (\ref{eq:prepro}) and (\ref{eq:proto})
where product of a transfer function and a scalar function of time is intended to
convey convolution.
%In the next section we introduce some preliminaries required to develop a
%learning control law inspired by the above equation.
%%%%%%%%%%%%%%%%%%%%%%%%%%%%%%%%%%%%%%%%%%%%%%%%%%%%%%%%%%%%%%%%%%%%%%%%%%%%%%%%
\subsection{Prototype Iterative Learning Control Algorithm}
\label{sec:proto}
\noindent
Analogous to (\ref{eq:filt}) we have in the ``lifted'' domain
\begin{eqnarray}
\label{eq:mat}
\mathbf G &=& \mathbf G^+\mathbf G^-
\end{eqnarray}
Note that the invertible part of the system $G^+(z^{-1})$ is represented by
the lower triangular matrix $\mathbf G^+$ consiting of its Markov parameters
(\ref{eq:mark}).
The all zero part of the plant $G^-(z^{-1})$ is
represented by the lower triangular banded Toeplitz matrix
\begin{eqnarray}
\label{eq:nu}
\mathbf G^-_{ij} &=& g_{i-j}, \quad 0\leq i-j\leq nu \nonumber \\
&=& 0, \quad otherwise
\end{eqnarray}
Hence the plant output defined in (\ref{eq:plant})
\[ \mathbf y_k = \mathbf G\mathbf u_k =\mathbf G^-\mathbf G^+\mathbf u_k =
\mathbf G^-\mathbf u'_k \]
where we have introduced the variable $\mathbf u' = \mathbf G^+ \mathbf u$.
Analogous to (\ref{eq:proto}) we propose the prototype learning law
\begin{eqnarray}
\label{eq:protoilc}
\mathbf u'_{k+1} = \mathbf u'_k + \alpha(\mathbf G^-)^T\mathbf e_k
\end{eqnarray}
We will address the stability of the prototype ILC in section~\ref{sec:stab}.
Assuming it is indeed stable the fixed point of iterations is given by
\begin{eqnarray}
(\mathbf G^-)^T\mathbf e_\infty &=& 0 \nonumber \\
(\mathbf G^-)^T \mathbf G^- \mathbf u'_\infty &=& (\mathbf G^-)^T\mathbf r
\end{eqnarray}
Hence we see that the prototype ILC, if stable, will converge to
the optimal solution of the least squares minimization problem
\begin{eqnarray}
\min_{\mathbf u'} ||\mathbf r - \mathbf G^-\mathbf u'||_2
\end{eqnarray}
\noindent

%%%%%%%%%%%%%%%%%%%%%%%%%%%%%%%%%%%%%%%%%%%%%%%%%%%%%%%%%%%%%%%%%%%%%%%%%%%%%%%%
\section{Modified Repetitive Control}
\label{sec:modrep}
\noindent
To incorporate robustness in the presence of plant model uncertainty,
the prototype repetitive control (\ref{eq:proto}) was
modified as follows \cite{Tom89}
\begin{eqnarray}
\label{eq:mod}
u'(t) &=& \left[z^{-nq}Q(z,z^{-1})\right]\{u'(t-N+nq) \nonumber\\
      &+& \alpha \left[z^{-nu}G^-(z)\right]e(t-N+d+nu+nq)\}
\end{eqnarray}
where $Q(z,z^{-1})$ is a zero-phase low pass filter(\ref{eq:qfil}).
The corresponding matrix in the ``lifted'' domain is
\begin{eqnarray}
\label{eq:nq}
\mathbf Q_{ij} &=& q_{|i-j|}, \quad |i-j|\leq nq \nonumber \\
&=& 0, \quad otherwise
\end{eqnarray}
We wish to establish a learning control law inspired by (\ref{eq:mod}) which would
result in a easy to check stability condition. In the next section we introduce the
central focus of this paper viz., the modified repetitive ILC algorithm. We have taken
great care in setting up the proposed algorithm so as to result in an elegant and easy to
check condition for iteration convergence.

%%%%%%%%%%%%%%%%%%%%%%%%%%%%%%%%%%%%%%%%%%%%%%%%%%%%%%%%%%%%%%%%%%%%%%%%%%%%%%%%
\subsection{Modified Repetitive ILC Algorithm}
\label{sec:blt}
\noindent
First we define learning gain filters $Q_u$ and
$Q_e$, of the form (\ref{eq:qfil}),
of orders $nq^u$ and $nq^e$ respectively.
We redefine the \textit{extended} plant input, output and reference signals by
\[
\begin{array}{ccl}
\mathbf u_k &=& [\begin{array}{ccc}u(0)& \ldots &
                           u(n+2nu-1)\end{array}]^T \\
\mathbf y_k &=& [\begin{array}{ccc}y(1)& \ldots &
                           y(n+2nu)\end{array}]^T \\
\mathbf r &=& [\begin{array}{ccc}r(1)& \ldots &
                           r(n+2nu)\end{array}]^T
\end{array}
\]
where we have assumed the relative degree $d=1$ without any loss of generality.
Now we zero-pad the control input as follows
\[
\begin{array}{ccll}
\mathbf u'_k(t) &=& 0,& t=0\ldots nu-1 \\
&=& \bar{\mathbf u}_{k}(t-nu),& t=nu\ldots n+nu-1 \\
&=& 0,& t=n+nu\ldots n+2nu-1
\end{array}
\]
which gives us the modified plant model
\[ \mathbf y_k = \mathbf G^-\mathbf u'_k = \mathbf G^-\mathbf N \bar{\mathbf u}_k\]
where the ``zero-padding'' matrix $\mathbf N$ is defined according to
\[\mathbf N_{n+2nu,n} = \left[\begin{array}{c}\mathbf 0_{nu,n} \\
                                      \mathbf I_{n,n} \\
                                      \mathbf 0_{nu,n}
                             \end{array} \right]\]
%where we have split the lower triangular matrix $\mathbf G^-$ into three blocks
%$\mathbf G^-_{n+2c,n+2c} = \left[\triangle_{n+2c,c} \quad | \bar{\mathbf G}_{n+2c,n}
%\quad | \nabla_{n+2c,c}\right]$.
Now we propose the learning control law
\begin{eqnarray}
\label{eq:law}
\bar{\mathbf u}_{k+1} = \mathbf Q_u\bar{\mathbf u}_k + \mathbf F\mathbf e_k
\end{eqnarray}
where $\mathbf Q_u$ is the $n\times n$ matrix representation of $Q_u$ as in (\ref{eq:nq}).
We can not derive an error propagation relation
because $\mathbf G^-\mathbf N$ and $\mathbf Q_u$ do not commute in general.
Instead we shall use a state space approach to establish stability where $\bar{\mathbf u}$,
$\mathbf r$ and $\mathbf e$ will
serve as the system state, input and output respectively.
%This general update allows for \textit{anticausal} (in time) filtering of the control and error
%signals from the previous cycle. Since we have the entire trial record $\bar u_k$
%and $e_k$ the update law can be higher order in time.
%In \cite{Tsao04} we start with the restricted choice $Q_u= I_{n,n}$
%and derive the required sufficiency conditions using a different approach.
%% Quote Moore paper for choice of Q and conexion to fixed point
%%%%%%%%%%%%%%%%%%%%%%%%%%%%%%%%%%%%%%%%%%%%%%%%%%%%%%%%%%%%%%%%
From (\ref{eq:errdef}) and (\ref{eq:law}) we have the control propagation equation
\begin{eqnarray}
\label{conprop}
%\begin{array}{ccl}
\bar{\mathbf u}_{k+1} &=& \mathbf Q_u\bar{\mathbf u}_k + \mathbf F\mathbf e_k \nonumber\\
&=& \mathbf Q_u\bar{\mathbf u}_k + \mathbf F(\mathbf r-\mathbf y_k) \nonumber\\
&=& \mathbf Q_u\bar{\mathbf u}_k - \mathbf F\mathbf G^-\mathbf u'_k + \mathbf F\mathbf r \nonumber\\
&=& (\mathbf Q_u-\mathbf F\mathbf G^-\mathbf N)\bar{\mathbf u}_k + \mathbf F\mathbf r
\end{eqnarray}
Now we choose $\mathbf F=\alpha \mathbf N^T(\mathbf G^-)^T \mathbf Q_e$
where $\mathbf Q_e$ is the $(n+2nu)\times (n+2nu)$ matrix representation of $Q_e$.
we end up with the state transition matrix
$\mathbf A = \mathbf Q_u-\alpha \mathbf N^T(\mathbf G^-)^T\mathbf Q_e\mathbf G^-\mathbf N$.
As was observed in \cite{Kev89} if $\mathbf Q_u\neq \mathbf I_{n,n}$ the
error vector no longer converges to the zero vector.
A necessary and sufficient condition for iteration convergence is
that the state transition matrix $\mathbf A$ be stable.
It is worth noting that the ``zero-padding'' matrix
$\mathbf N$ and the learning gains $\mathbf Q_u$ and $\mathbf Q_e$ have all been chosen carefully such that $\mathbf A $ has a SBT structure. It is this special structure that establishes the vital link between the ``lifted'' and frequency domain stability conditions.

%%%%%%%%%%%%%%%%%%%%%%%%%%%%%%%%%%%%%%%%%%%%%%%%%%%%%%%%%%%%%%%%%%%%%%%%%%%%%%%%%
\section{Stability Analysis}
\label{sec:stab}
\subsection{Symmetric Banded Toeplitz Structure}
\noindent
It can be shown that the state transition matrix $\mathbf A$
has the SBT structure below (c.f. Appendix~\ref{sec:sbltproof})
\[
\left[\begin{array}{ccccccccccc}
      a_0 & a_1 & \ldots & a_r&&&&&&&\\
      a_1 & a_0 &&&&&&&&&\\
      %\vdots &&&&&&0&&&&\\
      \vdots&&\ddots &&& \ddots &&&0&&\\
      a_r &&&&&&&&&&\\
      & \ddots &&&&&&&&&\\
      && a_r & \ldots & a_1 & a_0 & a_1 & \ldots & a_r &&\\
      &&&&&&&&&&\\
      &&0&&&&&\ddots&&&a_r\\
      &&&&&&&&&&\vdots\\
      &&&&&&&&&a_0&a_1\\
      &&&&&&&a_r&\ldots&a_1 & a_0
\end{array}\right]
\]
where $r=max(nu,nq^e+nu)$.
The entries of the matrix above, $a_i$, have been explicitly derived in Appendix~\ref{sec:Acoeff}.
%%%%%%%%%%%%%%%%%%%%%%%%%%%%%%%%%%%%%%%%%%%%%%%%%%%%%%%%%%%%%%%%%%%%%%%%%%%%%%%%
\subsection{Connections to Frequency Domain Stability Condition}
\label{sec:conn}
\noindent
%\[e_{k+1} = Ae_k\]
A necessary and sufficient condition for iteration convergence is
\begin{eqnarray}
\max_m|\lambda_m(\mathbf A)| < 1, \quad m=1\ldots n
\end{eqnarray}
where $\lambda_m$ stands for the $m^{th}$ eigenvalue of $\mathbf A$. This is the true stability condition (using the terminology in \cite{Long00}). Now let the discrete time noncausal filter representation of $\mathbf A$ be
%\[\mathbf{E}_{k+1}(z^{-1})=\mathbf{A}(z,z^{-1})\mathbf{E}_k(z^{-1})\]
%where the discrete time non-causal filter
\[A(z,z^{-1})=a_0+\sum_{k=1}^r a_k(z^k+z^{-k})\]
where the discrete time transform variable is evaluated at $z=e^{-j\omega T}$. Now
for sufficiency we look at the $H_\infty$ norm condition \cite{Owe92,Owe96}
\[
|A(e^{-j\omega T})|<1 \quad \forall \quad \omega T \in [0,\pi]\]
which translates to
\begin{eqnarray}
\label{eq:hinfnorm}
|a_0 + 2\sum_{k=1}^r a_k\cos(k\theta)|<1 \quad \forall \quad \theta \in [0,\pi]
\end{eqnarray}
We use a unique property of symmetric Toeplitz matrices (c.f. Section 4.2, Lemma 4.1 \cite{Gray}) viz.,
\begin{eqnarray}
\max_m|\lambda_m(\mathbf A)|<\max_{\theta\in [0,\pi]}|a_0 + 2\sum_{k=1}^r a_k\cos(k\theta)|
\end{eqnarray}
which gives us the required \emph{key} link between frequency domain and the ``lifted'' domain. For the choice of learning matrices we made, the above translates to the easy to check $H_\infty$ norm condition
\begin{eqnarray}
\label{eq:stab}
\left|Q_u(z,z^{-1})-\alpha Q_e(z,z^{-1})G^-(z^{-1})
G^-(z)\right|_\infty &<& 1
\end{eqnarray}
Using the terminology in \cite{Long00}, the above would be the approximate stability condition.
For the special choice $Q_u = Q_e = 1$ we get back the prototype ILC
(c.f. section~\ref{sec:proto})
with the simple stability condition
\begin{eqnarray}
\label{eq:stabsim}
\left|1-\alpha G^-(z^{-1})G^-(z)\right|_\infty &<& 1
\end{eqnarray}
which holds iff $\alpha \in [0,2]$ \cite{Tom90}.
We immediately note that for a general system it is easier to satisfy (\ref{eq:stab}) because
of the extra degrees of freedom in choosing $Q_u$ and $Q_e$.
%Also when
%$Q_u = Q_e$ the stability condition is identical to the robust stability
%established in \cite{Tom89}.
%We infer from this analysis that the learning filters $\mathbf Q_u$ and $\mathbf Q_e$
%are be chosen so as to satisfy (\ref{eq:modstab}) which in turn gaurantees iteration convergence
%in the ``lifted'' domain.
\\\indent
The above result is only a sufficient condition for error convergence.
The main obstacle in proving necessity is the absence of explicit formulae for the
eigenvalues of the matrix $\mathbf A$. To show that the condition is not overly
conservative we shall first perturb $\mathbf A$ to convert it into a
circulant matrix \cite{Gray} and second look at the special case when $\mathbf A$ is
a tri-diagonal matrix. In both situations we have the luxury of being able to compute the eigenvalues explicitly and show that the condition is fairly stringent.
Not surprisingly, efficient computation of the eigenvalues of a SBT matrix is in itself a well researched topic \cite{Bini,Arb}.
%%%%%%%%%%%%%%%%%%%%%%%%%%%%%%%%%%%%%%%%%%%%%%%%%%%%%%%%%%%%%%%%%%%%%%%%%%%%%%%%
\subsection{Circulant Matrix Approximation to the Transition Matrix}
\noindent
We can approximate $\mathbf A$ by embedding it in the circulant matrix $\tilde{\mathbf A}$ below \cite{Arb}
\[
\left[\begin{array}{ccccccccccc}
      a_0 & a_1 & \ldots & a_r &&&&&a_r&\ldots&a_1\\
      a_1 & a_0 &&&&&&&&\ddots&\vdots\\
      %\vdots &&&&&&0&&&&\\
      \vdots&&\ddots &&& \ddots &&&0&&a_r\\
      a_r &&&&&&&&&&\\
      & \ddots &&&&&&&&&\\
      && a_r & \ldots & a_1 & a_0 & a_1 & \ldots & a_r &&\\
      &&&&&&&&&&\\
      &&0&&&&&\ddots&&&a_r\\
      a_r&&&&&&&&&&\vdots\\
      \vdots&\ddots&&&&&&&&a_0&a_1\\
      a_1&\ldots&a_r&&&&&a_r&\ldots&a_1 & a_0
\end{array}\right]
\]
We note this is not a bad approximation considering that in practice $n\gg r$ and the
added terms do not perturb the eigenvalues too much. In fact we have for $n\gg r$ \cite{Gray}
\[\max_m|\lambda_m(\mathbf A)|\simeq\max_m|\lambda_m(\tilde{\mathbf A})|\]
The eigenvalues of $\tilde{\mathbf A}$ are given by the Discrete Fourier Transform (DFT)
of the first row of the matrix \cite{Dim,Gray}
\[
\begin{array}{ccl}
\lambda_m(\tilde{\mathbf A})&=& a_0+\sum_{k=1}^r a_k\left(e^{-\frac{2\pi ik(m-1)}
   {n}} + e^{-\frac{2\pi i(n-k)(m-1)}{n}}\right) \\
& = & a_0 + 2\sum_{k=1}^r a_k\cos\left(\frac{2\pi k(m-1)}{n}\right),
    \quad m=1\ldots n
\end{array}
\]
Hence the condition for convergence becomes
\[\max_m|a_0 + 2\sum_{k=1}^r a_k\cos(k\theta_m)|<1,\quad
  \theta_m=\frac{2\pi (m-1)}{n}\]
When $\theta_m\in [0,\pi]$ we readily have the sufficient condition for iteration
convergence from (\ref{eq:hinfnorm}). When $\theta_m\in [\pi,2\pi]$
let $\theta_m=2\pi-\theta$ where $\theta\in [0,\pi]$. Then we have
\[\cos(k\theta_m)=\cos(k2\pi-k\theta)=\cos(k\theta)\]
which gives us the required sufficient condition for iteration convergence.
In the limiting case when
$n\to\infty$, we see that $\theta_m$ spans the entire range $(0,\pi)$ and
hence the norm condition (\ref{eq:hinfnorm}) becomes both necessary and sufficient.
%The condition becomes closer to being necessary for large $n$.

%since we can shoose $\theta$ such that
%$\cos(\theta)=\frac{sgn(a_1)}{sgn(a_0)}$
%When $Q\neq I$ it is not possithe authors believe that the
%condition for monotonic convergence (\ref{eq:mon2}) would have a similar form
%involving the terms in the $Q$ matrix.
%%%%%%%%%%%%%%%%%%%%%%%%%%%%%%%%%%%%%%%%%%%%%%%%%%%%%%%%%%%%%%%%%%%%%%%%%%%%%%%%%%%%%%%%%%
\subsection{Three Banded Transition Matrix}
\noindent
If the system were such that $r=1$ then the resulting state transition matrix is
tridiagonal with eigenvalues given by \cite{Smi}
\[\lambda_m(\mathbf A) = a_0 + 2a_1\cos(\theta_m), \quad \theta_m=\left(\frac{m\pi}{n+1}\right),
\quad m=1\ldots n\]
From (\ref{eq:hinfnorm}) we have the $H_\infty$ norm condition
\begin{eqnarray}
|a_0 + 2a_1\cos(\theta)|<1 \quad \forall \quad \theta \in [0,\pi]
\end{eqnarray}
which gives a sufficient condition for iteration convergence. In the limiting case when
$n\to\infty$, we again see that $\theta_m$ spans the entire range $(0,\pi)$ and hence the
norm condition becomes both necessary and sufficient.
%We have equivalence between
%monotonic convergence (\ref{eq:mon2}) and the norm bound (\ref{eq:hinfnorm}) since
%\begin{eqnarray}
%|a_0 + 2a_1\cos(\phi)| <1 \Leftrightarrow |a_0| + 2|a_1|<1
%\end{eqnarray}
%where we choose $\phi$ to be $0$($\pi$) when $a_0$ and $a_1$ have the same(different) sign.
%%%%%%%%%%%%%%%%%%%%%%%%%%%%%%%%%%%%%%%%%%%%%%%%%%%%%%%%%%%%%%%%%%%%%%%%%%%%%%%%%%%%%%%%%%
\section{Monotonic Error Convergence}
\label{sec:mon}
\noindent
For the special choice $\mathbf Q_u = \mathbf I_{n,n}$ the modified
Repetitive ILC algorithm (see section~\ref{sec:blt}) simplifies to
\begin{eqnarray}
\bar{\mathbf u}_{k+1} &=& \bar{\mathbf u}_k + \mathbf F\mathbf e_k \nonumber \\
                      &=& (\mathbf I - \mathbf F\mathbf G^-\mathbf N)\bar{\mathbf u}_k
                      + \mathbf F\mathbf r \nonumber \\
                      &=& \mathbf A\bar{\mathbf u}_k + \mathbf F\mathbf r
\end{eqnarray}
Now if assume the initial conditions $\bar{\mathbf u}_0 = \mathbf 0$ and $\mathbf e_0 = \mathbf r$
we can show by induction that
\begin{eqnarray}
\mathbf F\mathbf e_k = \mathbf A^k\mathbf F\mathbf r
\end{eqnarray}
which gives us the error propagation equation
\begin{eqnarray}
\mathbf F \mathbf e_{k+1} &=& \mathbf A\mathbf F \mathbf e_k
\end{eqnarray}
Since $\mathbf A$ is symmetric by definition we have,
\begin{eqnarray}
||\mathbf A||^2_2 &=& \lambda_{max}(\mathbf A^T\mathbf A) \nonumber \\
&=& \lambda_{max}(\mathbf A^2) \nonumber \\
&\leq& ||\mathbf A^2||_1 \nonumber \\
&\leq& ||\mathbf A||^2_1
\end{eqnarray}
where we have used the fact that the spectral radius of a matrix is
bounded by its 1-norm.
For monotonic convergence in the p-norm (p=1,2 or $\infty$) we need
$||\mathbf F\mathbf e_{k+1}||_p \leq ||\mathbf F\mathbf e_k||_p$ for all $k$.
\begin{eqnarray}
||\mathbf F\mathbf e_{k+1}||_p &=& ||\mathbf A\mathbf F\mathbf e_k||_p \nonumber \\
&\leq& ||\mathbf A||_p||\mathbf F\mathbf e_k||_p \nonumber \\
&\leq& ||\mathbf A||_1||\mathbf F\mathbf e_k||_p
\end{eqnarray}
where we have used the relation
$||\mathbf A||_2 \leq||\mathbf A||_\infty=||\mathbf A||_1$.
Since the induced-1 norm equals the maximal column sum we have the
sufficient condition for monotonic convergence
\begin{eqnarray}
\label{eq:mon2}
||\mathbf A||_1=|a_0| + 2\sum_{k=1}^{nq_e+nu}|a_k|< 1
\end{eqnarray}
Hence we conclude that the error converges to the fixed point of the
iteration $\mathbf F \mathbf e_\infty = \mathbf 0$ if $\mathbf A$ is stable. In
addition if (\ref{eq:mon2}) is satisfied it does so in a monotonic decreasing
fashion. As expected the convergence condition (\ref{eq:hinfnorm}) is implied by the
stronger monotonicity condition (\ref{eq:mon2}) as shown below.
\begin{eqnarray}
|a_0 + 2\sum_{k=1}^{nq_e+nu} a_k\cos(k\theta)| &\leq&
|a_0| + 2\sum_{k=1}^{nq_e+nu} |a_k\cos(k\theta)| \nonumber \\
&\leq& |a_0| + 2\sum_{k=1}^{nq_e+nu} |a_k||\cos(k\theta)| \nonumber \\
&\leq& |a_0| + 2\sum_{k=1}^{nq_e+nu} |a_k|
\end{eqnarray}
Also the error vector $\mathbf e_\infty = \mathbf 0$ only when the system is
fully invertible i.e., $nu=0$ and the resultant $\mathbf F$ is a nonsingular matrix.

%%%%%%%%%%%%%%%%%%%%%%%%%%%%%%%%%%%%%%%%%r%%%%%%%%%%%%%%%%%%%%%%%%%%%%%%%%%%%%%%%
\section{Numerical Example}
\label{sec:ex}
\begin{figure}[t]
\centerline{\scalebox{.55}{\includegraphics*[50mm,85mm][175mm,205mm]{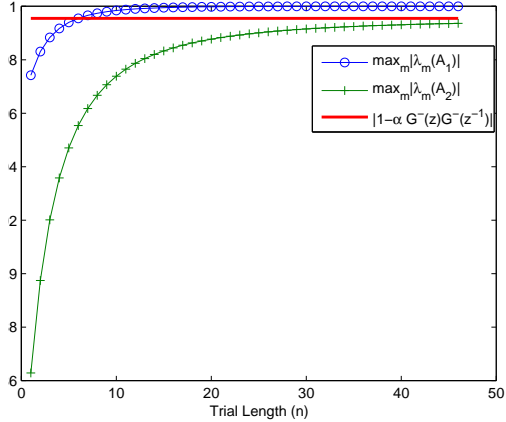}}}
\caption{Effect of ``zero-padding'' matrix on ILC stability}
\label{lamda}
\end{figure}

\noindent
We shall illustrate the reason behind the introduction of the ``zero-padding'' matrix $\mathbf N$ in section~\ref{sec:stab} with a simple example. Let us consider
the following non-minimum phase plant \cite{Kev01}
\[y(t+1) = -0.2y(t) + 0.0125y(t-1) + u(t) - 1.1u(t-1)\]
From the definitions earlier (c.f. Section \ref{sec:proto}) we have
$G^-(z^{-1}) = (g_0+g_1z^{-1})$ where $g_0=1$ and
$g_1=-1.1$. The number of unstable zeros $nu=1$. Now we use the proposed design methodology with
$Q_u(z,z^{-1})=Q_e(z,z^{-1})=1$ and $\alpha=0.45$. Hence we have the noncausal filter representation for the transition matrix given by
\begin{eqnarray*}
A(z,z^{-1}) &=& 1-\alpha G^-(z)G^-(z^{-1}) \\
&=& 1-\alpha\left(g_0^2+g_1^2-g_0g_1z-g_0g_1z^{-1}\right)
\end{eqnarray*}
This choice satisfies the convergence condition (\ref{eq:hinfnorm})
\begin{eqnarray}
\label{eq:bound1}
|a_0+2a_1cos(\theta)|<1,\quad \forall \theta\in[0,\pi]\nonumber\\
\Rightarrow\left|1-\alpha(g_0^2+g_1^2)\right| + 2|\alpha|\left|g_0g_1\right|=0.9955<1
\end{eqnarray}
since $a_0 = 1-\alpha\left(g_0^2+g_1^2\right)$ and $a_1 = -\alpha g_0g_1$.
We have the transition matrix without zero-padding
given by $\mathbf A_1 = \mathbf I - \alpha(\mathbf G^-)^T\mathbf G^-$.
Instead if we do the zero-padding we get the matrix
$\mathbf A_2 = \mathbf I-\alpha\mathbf N^T(\mathbf G^-)^T\mathbf G^-\mathbf N$. To understand the value of the zero-padding, one can easily check for trial length $n=3$,
\begin{eqnarray*}
\mathbf A_1 &=& \left[\begin{array}{ccc} 0.0055 &   0.4950 &        0\\
    0.4950 &   0.0055   & 0.4950\\
         0  &  0.4950  &  0.5500\end{array}\right]\\
\mathbf A_2 &=& \left[\begin{array}{ccc} 0.0055 &   0.4950 &        0\\
    0.4950 &   0.0055   & 0.4950\\
         0  &  0.4950  &  0.0055\end{array}\right]
         \end{eqnarray*}
The ``zero-padding'' matrix $\mathbf N$ ensures that the ($n,n$) entry of the
transition matrix is $0.0055$ thereby marking it a banded toeplitz matrix. Remarkable as it may seem, difference in this single entry results in $\lim_{n\rightarrow\infty}\max_m|\lambda_m(\mathbf A_1)|=1$ whereas $\lim_{n\rightarrow\infty}\max_m|\lambda_m(\mathbf A_2)|=0.9955$ i.e., the upper bound in (\ref{eq:bound1}) (see fig.~\ref{lamda}).
We conclude that the ``zero-padding'' matrix $\mathbf N$ is critical in ensuring stability of the learning algorithm. An injudicious choice of the learning matrices would result in instability despite the approximate stability condition (\ref{eq:hinfnorm}) being met!

%%%%%%%%%%%%%%%%%%%%%%%%%%%%%%%%%%%%%%%%%%%%%%%%%%%%%%%%%%%%%%%%%%%%%%%%%%%%%%%%

\section{Experimental and Simulation Results}
\begin{figure}[t]
\centerline{\scalebox{.425}{\includegraphics*[15mm,60mm][195mm,220mm]{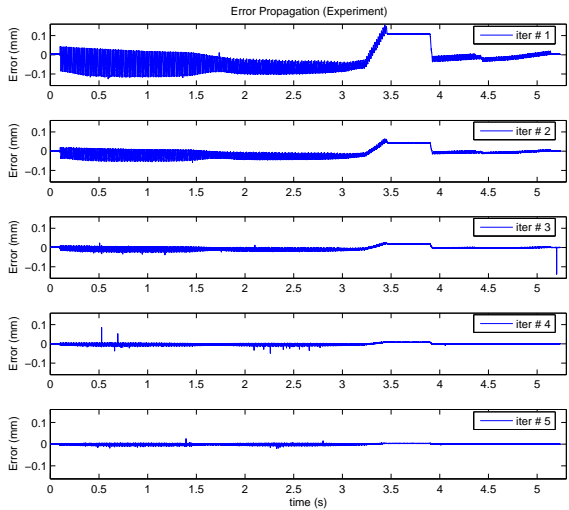}}}
\caption{Error Propagation for First Five Iterations (Experiment)}
\label{err}
\end{figure}

\begin{figure}[t]
\centerline{\scalebox{.525}{\includegraphics*[35mm,80mm][180mm,195mm]{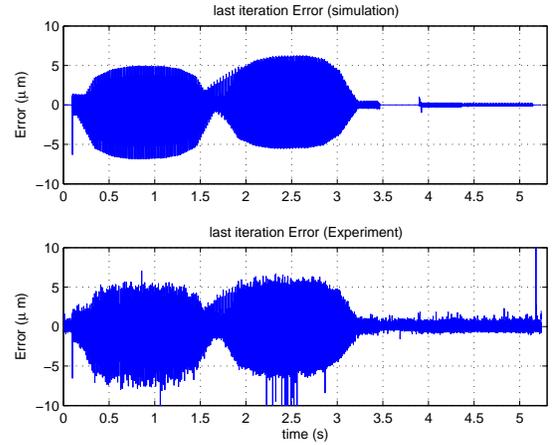}}}
\caption{Final Iteration Error Simulation (top) and Experiment (bottom)}
\label{final}
\end{figure}

%\subsection{Plant Model and Control Design}
\label{sec:exp}
\noindent
The plant to be controlled is a DC motor driven cutting tool with a stabilizing PD control.
The
system was identified in the discrete domain at a sampling frequency of T=15kHz.
\begin{figure}[t]
\centerline{\scalebox{.55}{\includegraphics*[30mm,80mm][190mm,200mm]{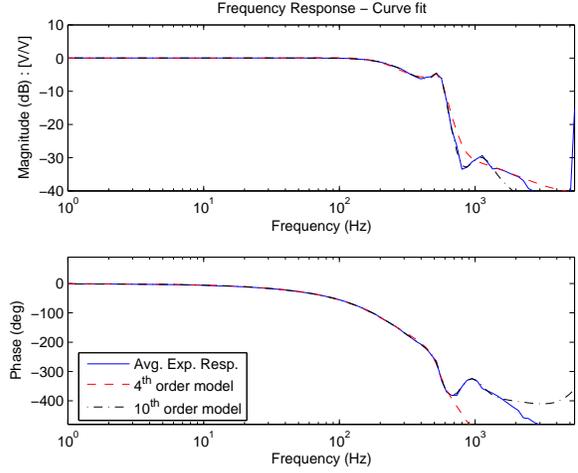}}}
\caption{PD Plant Experimental and Model Frequency Response}
\label{sysmodel}
\end{figure}
Figure~\ref{sysmodel} shows the PD plant and the curve fit model frequency response. We
have a $4^{th}$ order model which is used for control design and also a $10^{th}$ order model
which will be used in a later section. We will design a Iterative learning control to better
the
performance of the inner loop controller. As can be seen from the model frequency response
(Fig.~\ref{sysmodel}) this
system has a low bandwidth since the phase drops quite sharply ($100^o$ by 200 Hz).
%%%%%%%%%%%%%%%%%%%%%%%%%%%%%%%%%%%%%%%%%%%%%%%%%%%%%%%%%%%%%%%%%%%%%%%%%%%%%%%%%%%%
\begin{figure}[t]
\centerline{\scalebox{.45}{\includegraphics*[20mm,65mm][195mm,215mm]{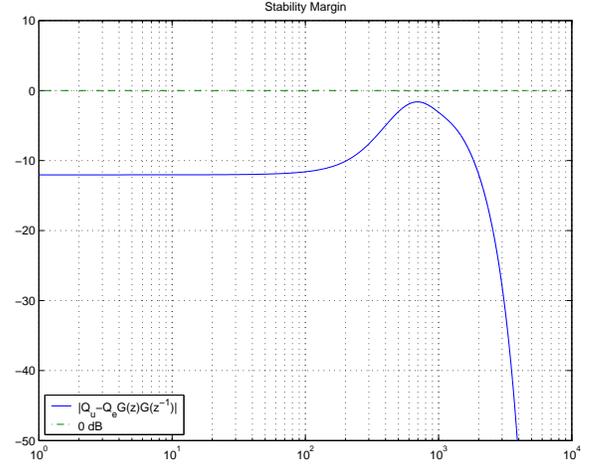}}}
\caption{Stability Condition $\left|Q_u(z,z^{-1})-\alpha Q_e(z,z^{-1})G^-(z^{-1})
G^-(z)\right|_\infty < 1$}
\label{stab}
\end{figure}
\begin{figure}[t]
\centerline{\scalebox{.45}{\includegraphics*[15mm,65mm][195mm,205mm]{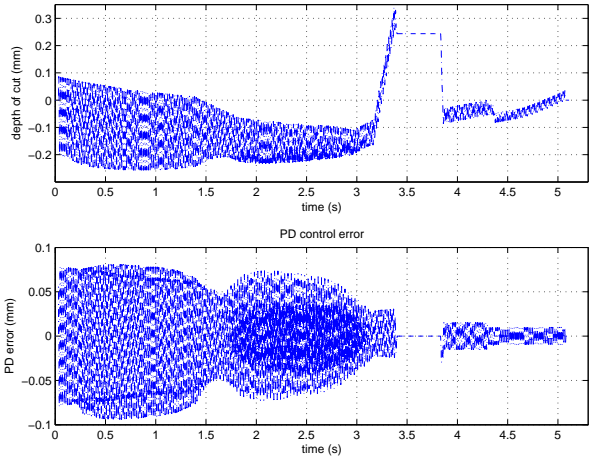}}}
\caption{Complex Tool Profile (top) with PD Tracking Error (bottom)}
\label{refcmd}
\end{figure}
\noindent
The top plot of Fig.~\ref{refcmd} shows the complex profile as a function of time.
What is
shown is the depth of cut that the tool tip has to achieve. The command profile has periodic
content
between $t=1s$ and $t=4s$
reflecting oval cross-section and flat regions between $t=3.5s$ and $t=3.8s$
reflecting round cross-section.
The profile is complex in the sense that although it looks periodic, the magnitude and phase
changes gradually making it
hard to follow. Also phase compensation \textit{apriori} based on the plant's frequency
response
(Fig.~\ref{sysmodel})
is not possible due to the above reason.
The bottom plot of Fig.~\ref{refcmd} shows the error for the PD plant.
The PD plant is not able to track this profile due to significant phase lag in the
frequencies contained in the
reference. We choose learning gain $\alpha=0.75$ and zero phase low pass filters
$Q_u$, $Q_e$ of orders 16 and 32 respectively.
Fig.~\ref{stab} shows that the stability condition is met for this choice of learning matrices.
We start the experiment with $\mathbf u'_0=0$, $\mathbf e_0 = r$.
Figure~\ref{err} shows the tracking error for the first
five iterations. We see that the error decreases uniformly over the entire time axis.
Fig.~\ref{final} shows the $10^{th}$ iteration tracking error achieved by the scheme.
We also did a simulation with no model mismatch for comparison with the experimental
results . The error is within an acceptable
$\pm 10\mu m$ compared to $\pm 100\mu m$ achieved
with the PD controlled plant (see bottom plot of Fig.~\ref{refcmd}).

%%%%%%%%%%%%%%%%%%%%%%%%%%%%%%%%%%%%%%%%%%%%%%%%%%%%%%%%%%%%%%%%%%%%%%%%%%%%%%%%
\subsection{Improvement over ZPETC Feedforward Control}
\begin{figure}[t]
\centerline{\scalebox{.525}{\includegraphics*[30mm,80mm][195mm,195mm]{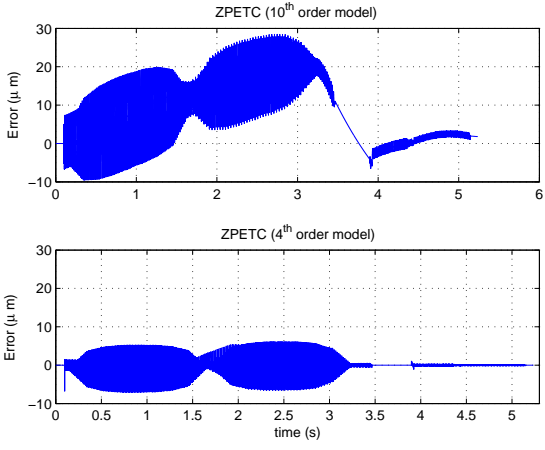}}}
\caption{$10^{th}$ order model (top) $4^{th}$ order model (bottom) ZPETC error}
\label{zpetc}
\end{figure}

\begin{figure}[t]
\centerline{\scalebox{.5}{\includegraphics*[30mm,75mm][195mm,215mm]{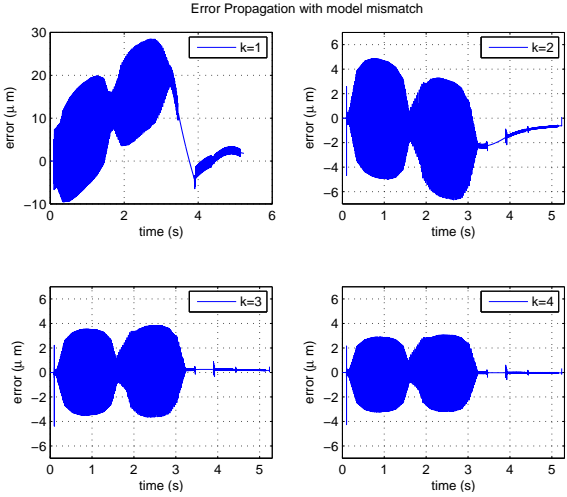}}}
\caption{Error Propagation with model mismatch (Simulation)}
\label{err1}
\end{figure}

\noindent
\label{sec:zpetc}
Assuming the system was at rest $\mathbf e_0=\mathbf r$ and $\mathbf u'_0=0$ we have
\[\mathbf u'_1=\alpha(\mathbf G^-)^T\mathbf r \Rightarrow
\mathbf u_1 = \alpha(\mathbf G^+)^{-1}(\mathbf G^-)^T\mathbf r\]
For the special choice of learning gain $\alpha = 1$
this is equivalent to the celebrated zero phase error tracking control (ZPETC) \cite{Tom85}.
Hence we infer that the first iteration of the prototype ILC (see section \ref{sec:proto})
corresponds to
stable inversion based feedforward compensation.
It is well known that feedforward control such as the ZPETC is sensitive to plant model
uncertainties.
But the prototype ILC has an inherent cycle-to-cycle feedback which makes it ``robust''
in a limited sense. As reported in \cite{Long02} zero phase ILC can be seen as repeated
application of the feedforward control leading to better performance despite the presence
of model uncertainties. Also the learning gain $\alpha$ is seldom set to the aggressive
value of $1$. Instead we use a lower gain and multiple iterations to derive the same effect
as feedforward control in a ideal plant model situation.
To verify the above claim, we did a simulation using a $4^{th}$ order model of the plant for
the control design. Then we applied ZPETC feedforward control to a more accurate $10^{th}$
order model (see fig.~\ref{sysmodel}). Figure~\ref{zpetc} shows that the feedforward
error is within $\pm 30 \mu m$ when there is a model mismatch (top) as against $\pm 5 \mu m$
when there is perfect model match (bottom).
Figure~\ref{err1} shows that we can improve upon the model mismatch
error using the proposed modified repetitive ILC and bring it within the
acceptable $\pm 5 \mu m$ region
in as few as four iterations.
%%%%%%%%%%%%%%%%%%%%%%%%%%%%%%%%%%%%%%%%%%%%%%%%%%%%%%%%%%%%%%%%%%%%%%%%%%%%%%%%%%%%%%%%%%%%%%%
\section{Conclusions}
\label{sec:con}
\noindent
%From the numerical simulation we see that there exists a tradeoff between achievable
%control and residual error based on the choice of the $Q$ filter in the proposed design.
We have introduced a noncausal ILC design based on zero-phase error compensation and filtering used in repetitive control design. The learning gains (matrices) are chosen carefully such that the resulting state transition matrix has a \textit{Symmetric Banded Toeplitz} (SBT) structure. This special structure has been exploited
to arrive at sufficient conditions for tracking error convergence in the frequency domain. Thus the ILC
matrix design is translated to repetitive control loop shaping filter design in the frequency domain. For a special choice of learning matrix we have also derived a stronger monotonic convergence condition. Furthermore we have shown that under plant model mismatch the proposed scheme improves upon stable inversion based feedforward control. The proposed methodology has been successfully demonstrated on a fast tool servo system used in precision machining applications.
%\\\indent
%Future
%work will focus on robust stability of the proposed algorithm in the presence
%of unmodeled plant dynamics.
%%%%%%%%%%%%%%%%%%%%%%%%%%%%%%%%%%%%%%%%%%%%%%%%%%%%%%%%%%%%%%%%%%%%%%%%%%%%%%%%
\appendices
\section{}
\label{sec:sbltproof}
\newtheorem{theorem}{Theorem}
\begin{theorem}{
The matrix $\mathbf A = \mathbf Q_u-\alpha \mathbf N^T(\mathbf G^-)^T\mathbf Q_e\mathbf G^-\mathbf N$ has a symmetric banded Toeplitz (SBT) structure}
\begin{IEEEproof}
Recall the definitions
\begin{eqnarray}
\mathbf G^-_{i,j} &=& g_{i-j}, \quad 0\leq i-j\leq nu \nonumber\\
&=& 0, \quad otherwise \nonumber \\
(\mathbf Q_e)_{i,j} &=& q^e_{|i-j|}, \quad |i-j|\leq nq \nonumber\\
&=& 0, \quad otherwise \nonumber \\
\mathbf N_{i,j} &=& 1, \quad 1\leq j=i-nu\leq n \nonumber \\
&=& 0, \quad otherwise \nonumber
\end{eqnarray}
Hence we have the coefficients of $(\mathbf Q_e\mathbf G^-)_{n+2nu,n+2nu}$ given by
\begin{eqnarray}
%(\mathbf {Q^eG}^-)_{i,j} &=& \sum_{k=max(i-nq,j)}^{min(i+nq,j+nu,n+2nu)}q_{|i-k|}g_{i-j}, \nonumber \\
b_{i,j} &=& \sum_{k=max(i-nq,j)}^{min(i+nq^e,j+nu,n+2nu)}q^e_{|i-k|}g_{i-j}, \nonumber \\
&&\quad\quad -nq^e\leq i-j\leq nq^e+nu \nonumber\\
&=& 0, \quad otherwise \nonumber
\end{eqnarray}
The coefficients of $(\mathbf Q_e\mathbf G^-\mathbf N)_{n+2nu,n}$ are given by
\begin{eqnarray}
c_{i,j} &=& \nonumber\\
b_{i,j+nu}
&=&\sum_{k=max(i-nq^e,j+nu)}^{min(i+nq^e,j+2nu)}q^e_{|i-k|}g_{k-j-nu}, \nonumber \\
&&\quad\quad nu-nq^e\leq i-j\leq nq^e+2nu \nonumber\\
&=& 0, \quad otherwise \nonumber
\end{eqnarray}
where we have used $j\leq n\Rightarrow j+2nu\leq n+2nu$.\\
The coefficients of $((\mathbf G^-)^T\mathbf Q_e\mathbf G^-\mathbf N)_{n+2nu,n}$ are given by
\begin{eqnarray}
d_{l,j} &=&
\sum_{i=1}^{n+2nu}(\mathbf G^-)_{i,l}c_{i,j} \nonumber \\
&=&\sum_{i=max(l,j+nu-nq^e)}^{min(l+nu,j+2nu+nq^e,n+2nu)}\{g_{i-l} \nonumber \\
&&\sum_{k=max(i-nq^e,j+nu)}^{min(i+nq^e,j+2nu)}q^e_{|i-k|}g_{k-j-nu}\}, \nonumber \\
&&\quad\quad -nq^e\leq l-j\leq nq^e+2nu \nonumber\\
&=& 0, \quad otherwise \nonumber
\end{eqnarray}
The coefficients of $(\mathbf N^T(\mathbf G^-)^T\mathbf Q_e\mathbf G^-\mathbf N)_{n,n}$ are given by
\begin{eqnarray}
e_{l,j} &=& \nonumber \\
d_{l+nu,j}
&=&\sum_{i=max(l+nu,j+nu-nq^e)}^{min(l+2nu,j+2nu+nq^e)}\{g_{i-l-nu} \nonumber \\
&&\sum_{k=max(i-nq^e,j+nu)}^{min(i+nq^e,j+2nu)}q^e_{|i-k|}g_{k-j-nu}\}, \nonumber \\
&&\quad\quad |l-j|\leq nq^e+nu \nonumber\\
&=& 0, \quad otherwise \nonumber
\end{eqnarray}
where we have used $l\leq n\Rightarrow l+2nu\leq n+2nu$.\\
We do a change of variables $i'=i-nu$ and $k'=k-nu$ to get
\begin{eqnarray}
e_{l,j}
&=&\sum_{i'=max(l,j-nq^e)}^{min(l+nu,j+nu+nq^e)}\{g_{i'-l} \nonumber \\
&&\sum_{k'=max(i'-nq^e,j)}^{min(i'+nq^e,j+nu)}q^e_{|i'-k'|}g_{k'-j}\}, \nonumber \\
&&\quad\quad |l-j|\leq nq^e+nu \nonumber\\
&=& 0, \quad otherwise \nonumber
\end{eqnarray}
We readily see that
\begin{eqnarray}
e_{l+1,j+1}
&=&\sum_{i'=max(l+1,j+1-nq^e)}^{min(l+1+nu,j+1+nu+nq^e)}\{g_{i'-l-1} \nonumber \\
&&\sum_{k'=max(i'-nq^e,j+1)}^{min(i'+nq^e,j+1+nu)}q^e_{|i'-k'|}g_{k'-j-1}\} \nonumber \\
&=&\sum_{i=max(l,j-nq^e)}^{min(l+nu,j+nu+nq^e)}\{g_{i-l} \nonumber \\
&&\sum_{k=max(i-nq^e,j)}^{min(i+nq^e,j+nu)}q^e_{|i-k|}g_{k-j}\} \nonumber \\
&=&e_{l,j} \nonumber
\end{eqnarray}
by doing a change of variable $i=i'-1$ and $k=k'-1$.
Hence $\mathbf N^T(\mathbf G^-)^T\mathbf Q_e\mathbf G^-\mathbf N$
(which is symmetric by definition) is a SBT matrix. Since $\mathbf Q_u$ is by definition a SBT matrix, it follows that $\mathbf A$ is also a SBT matrix.
\end{IEEEproof}
\end{theorem}
%%%%%%%%%%%%%%%%%%%%%%%%%%%%%%%%%%%%%%%%%%%%%%%%%%%%%%%%%%%%%%%%%%%%%%%%%%%%%

\section{}
\label{sec:Acoeff}
\noindent
\textbf{Entries of the state transition matrix A} -
\\\noindent
When $nq^u>nq^e+nu$ the
coefficients of $\mathbf A$ are given by
\begin{eqnarray}
f_{l,j}
&=&q^u_{|l-j|}-\alpha\sum_{i'=max(l,j-nq^e)}^{min(l+nu,j+nu+nq^e)}\{g_{i'-l} \nonumber \\
&&\sum_{k'=max(i'-nq^e,j)}^{min(i'+nq^e,j+nu)}q_{|i'-k'|}g_{k'-j}\}, \nonumber \\
&&\quad\quad |l-j|\leq nq^e+nu \nonumber\\
&=& q^u_{|l-j|},  \quad nq^e+nu<|l-j|\leq nq^u \nonumber\\
&=& 0, \quad otherwise \nonumber
\end{eqnarray}
When $nq^u\leq nq^e+nu$ the
coefficients of $\mathbf A$ are given by
\begin{eqnarray}
f_{l,j}
&=&q^u_{|l-j|}-\alpha\sum_{i'=max(l,j-nq^e)}^{min(l+nu,j+nu+nq^e)}\{g_{i'-l} \nonumber \\
&&\sum_{k'=max(i'-nq^e,j)}^{min(i'+nq^e,j+nu)}q_{|i'-k'|}g_{k'-j}\}, \nonumber \\
&&\quad\quad |l-j|\leq nq^u \nonumber\\
&=& -\alpha\sum_{i'=max(l,j-nq^e)}^{min(l+nu,j+nu+nq^e)}\{g_{i'-l} \nonumber \\
&&\sum_{k'=max(i'-nq^e,j)}^{min(i'+nq^e,j+nu)}q_{|i'-k'|}g_{k'-j}\}, \nonumber \\
&&\quad\quad nq^u<|l-j|\leq nq^e+nu \nonumber\\
&=& 0, \quad otherwise \nonumber
\end{eqnarray}
We get the entries in section~\ref{sec:stab} from the first row of the matrix
\[a_k = f_{1,k+1}, \quad 0\leq k\leq r\]
where $r=max(nq^u,nq^e+nu)$
%%%%%%%%%%%%%%%%%%%%%%%%%%%%%%%%%%%%%%%%%%%%%%%%%%%%%%%%%%%%%%%%%%%%%%%%%%%%%%%
\bibliographystyle{IEEEtran}
\bibliography{IEEEabrv,ref_ilc}

\begin{thebibliography}{10}
\providecommand{\url}[1]{#1}
\csname url@samestyle\endcsname
\providecommand{\newblock}{\relax}
\providecommand{\bibinfo}[2]{#2}
\providecommand{\BIBentrySTDinterwordspacing}{\spaceskip=0pt\relax}
\providecommand{\BIBentryALTinterwordstretchfactor}{4}
\providecommand{\BIBentryALTinterwordspacing}{\spaceskip=\fontdimen2\font plus
\BIBentryALTinterwordstretchfactor\fontdimen3\font minus
  \fontdimen4\font\relax}
\providecommand{\BIBforeignlanguage}[2]{{%
\expandafter\ifx\csname l@#1\endcsname\relax
\typeout{** WARNING: IEEEtran.bst: No hyphenation pattern has been}%
\typeout{** loaded for the language `#1'. Using the pattern for}%
\typeout{** the default language instead.}%
\else
\language=\csname l@#1\endcsname
\fi
#2}}
\providecommand{\BIBdecl}{\relax}
\BIBdecl

\bibitem{Ari}
S.~Arimoto, S.~Kawamura, and F.~Miyazaki, ``Bettering operation of robots by
  learning,'' \emph{J. of Robotic Systems}, vol.~1, no.~2, pp. 123--140, 1984.

\bibitem{Kevb}
K.~L. Moore, \emph{Iterative Learning Control for Deterministic Systems}, ser.
  Advances in Industrial Control.\hskip 1em plus 0.5em minus 0.4em\relax
  Springer-Verlag, 1993.

\bibitem{Bien}
Z.~Bien and J.-X. Xu, \emph{Iterative Learning Control - Analysis, Design,
  Integration and Applications}.\hskip 1em plus 0.5em minus 0.4em\relax Boston,
  MA: Kluwer Academic Publishers, 1998.

\bibitem{Kev98}
K.~L. Moore, ``Iterative learning control: A expository overview,''
  \emph{Applied and Computational Controls, Signal Processing and Circuits},
  vol.~1, no.~1, 1998.

\bibitem{Kev07}
H.-S. Ahn, Y.~Chen, and K.~L. Moore, ``Iterative learning control: Brief survey
  and categorization,'' \emph{IEEE Trans. systems, man and cybernetics-part c:
  applications and reviews}, vol.~37, no.~6, pp. 1099--1121, 2007.

\bibitem{Kev89}
K.~L. Moore, M.~Dahleh, and S.~Bhattacharyya, ``Iterative learning for
  trajectory control,'' in \emph{Proc. IEEE Conf. Decision and Control}, Tampa,
  FL, Dec 1989, pp. 860--865.

\bibitem{Owe92}
J.~B. Edwards and D.~H. Owens, \emph{Stability Analysis for Linear Repetitive
  Processes}, ser. Lecture Notes on Control and Information Sciences.\hskip 1em
  plus 0.5em minus 0.4em\relax NJ: Springer-Verlag, 1992.

\bibitem{Owe96}
N.~Amann, D.~H. Owens, E.~Rogers, and A.~Wahl, ``An $h_\infty$ approach to
  linear iterative learning control design,'' \emph{Int. J. Adaptive Control
  and Signal Processing}, vol.~10, pp. 767--781, 1996.

\bibitem{Dim}
D.~Gorinevsky, ``Loop-shaping for iterative learning control of batch
  processes,'' \emph{IEEE Control Systems Mag.}, vol.~22, no.~6, pp. 55--65,
  Dec 2002.

\bibitem{Nor02}
M.~M. Norl\"of and S.~Gunnarsson, ``Time and frequency domain convergence
  properties in iterative learning control,'' \emph{Int. J. Control}, vol.~75,
  no.~14, pp. 1114--1126, 2002.

\bibitem{Kev}
K.~L. Moore and Y.~Chen, ``On the monotonic convergence of high-order iterative
  learning updating laws,'' in \emph{Proc. IFAC Congress}, Barcelona, Spain,
  July 2002.

\bibitem{Kev02}
J.~H\"at\"onen, K.~L. Moore, and D.~H. Owens, ``An algebraic approach to
  iterative learning control,'' in \emph{Proc. IEEE Int. Symposium on
  Intelligent Control}, Vancouver, Canada, 2002, pp. 37--42.

\bibitem{Kev03}
K.~L. Moore and Y.~Chen, ``A seperative high-order framework for monotonic
  convergent iterative learning controller design,'' in \emph{Proc. American
  Control Conf.}, Denver, CO, June 2003, pp. 3644--3649.

\bibitem{Nor}
M.~Norl\"of, ``Comparative study on first and second order ilc - frequency
  domain analysis and experiments,'' in \emph{Proc. IEEE Conf. Decision and
  Control}, Sydney, Australia, Dec 2000, pp. 3415--3420.

\bibitem{Jian}
J.-X. Xu and Y.~Tan, ``On the convergence speed of a class of higher-order ilc
  schemes,'' in \emph{Proc. IEEE Conf. Decision and Control}, Orlando, FL, Dec
  2001, pp. 4932--4937.

\bibitem{Tom85}
M.~Tomizuka, ``Zero phase error tracking algorithm for digital control,''
  \emph{J. Dynamic Systems, Measurement and Control}, vol. 109, no.~1, pp.
  87--92, Mar 1985.

\bibitem{Long00}
R.~W. Longman, ``Iterative learning control and repetitive control for
  engineering practise,'' \emph{Int. J. Control}, vol.~73, no.~10, pp.
  930--954, 2000.

\bibitem{Long02}
H.~Elci, R.~W. Longman, M.~Q. Phan, J.-N. Juang, and R.~Ugoletti, ``Simple
  learning control made practical by zero-phase filtering: Applications to
  robotics,'' \emph{IEEE Trans. Circuits and Systems}, vol.~49, no.~6, pp.
  753--767, June 2002.

\bibitem{Long03}
S.~Songschon and R.~W. Longman, ``Comparison of the stability boundary and the
  frequency response stability condition in learning and repetitive control,''
  \emph{Int. J. Appl. Math. Comput. Sci}, vol.~13, no.~2, pp. 169--177, 2003.

\bibitem{Xu06}
J.-X. Xu and R.~Yan, ``On repetitive learning control for periodic tracking
  tasks,'' \emph{IEEE Trans. Automatic Control}, vol.~51, no.~11, pp.
  1842--1848, 2006.

\bibitem{Okko}
D.~D. Roover, O.~H. Bosgra, and M.~Steinbuch, ``Internal-model-based design of
  repetitive and iterative learning controllers for linear multivariable
  systems,'' \emph{Int. J. Control}, vol.~73, no.~10, pp. 914--929, 2000.

\bibitem{Chen}
Y.~Chen, K.~L. Moore, J.~Yu, and T.~Zhang, ``Iterative learning control and
  repetitive control in hard disk drive industry - a tutorial,'' in \emph{Proc.
  IEEE Conf. Decision and Control}, San Diego, CA, Dec 2006, pp. 1--15.

\bibitem{Long94a}
H.~Elci, R.~W. Longman, M.~Q. Phan, J.-N. Juang, and R.~Ugoletti, ``Automated
  learning control through model updating for precision motion control,''
  \emph{Adaptive Structures and Composite Materials: Analysis and Application},
  vol.~45, pp. 299--314, 1994.

\bibitem{Owe98}
N.~Amann, D.~Owens, and E.~Rogers, ``Predictive optimal iterative learning
  control,'' \emph{Int. J. Control}, vol.~69, no.~2, pp. 203--226, 1998.

\bibitem{Verw}
M.~H. Verwoerd, G.~Meinsma, and T.~J. de~Vries, ``On the use of noncausal lti
  operators in iterative learning control,'' in \emph{Proc. IEEE Conf. Decision
  and Control}, Las Vegas, NV, Dec 2002, pp. 3362--3366.

\bibitem{Verw05}
\BIBentryALTinterwordspacing
M.~H. Verwoerd, ``Iterative learning control - a critical review,'' Ph.D.
  dissertation, University of Twente, Netherlands, 2005. [Online]. Available:
  \url{http://www.ce.utwente.nl/rtweb/publications/2005/pdf-files/101CE2005\_Verwoerd.pdf}
\BIBentrySTDinterwordspacing

\bibitem{Wu09}
Y.~Wu, Q.~Zou, and C.~Su, ``A current cycle feedback iterative learning control
  approach for afm imaging,'' \emph{IEEE Trans. Nanotechnology}, vol.~8, no.~4,
  pp. 515--527, 2009.

\bibitem{Leang}
K.~K. Leang and S.~Devasia, ``Iterative feedforward compensation of hysteresis
  in piezo positioners,'' in \emph{Proc. IEEE Conf. Decision and Control},
  Maui, Hawaii, Dec 2003, pp. 2626--2631.

\bibitem{Deva05}
S.~Tien, Q.~Zou, and S.~Devasia, ``Iterative control of dynamics-coupling
  effects in piezo-actuator for high-speed afm operation,'' \emph{IEEE Trans.
  Contr. Syst. Technol.,}, vol.~13, no.~6, pp. 921--931, 2005.

\bibitem{Kri3}
K.~Krishnamoorthy, C.~Y. Lin, and T.-C. Tsao, ``Design and control of a dual
  stage fast tool servo for precision machining,'' in \emph{Proc. IEEE Conf.
  Control Applications}, Taipei, Taiwan, 2004, pp. 742--747.

\bibitem{Long88}
M.~Phan and R.~Longman, ``A mathematical theory of learning control for linear
  discrete multivariable systems,'' in \emph{Proc. AIAA/AAS Astrodynamics
  Conf.}, Minneapolis, MN, Aug 1988, pp. 740--746.

\bibitem{Gray}
\BIBentryALTinterwordspacing
R.~M. Gray, ``Toeplitz and circulant matrices: A review,'' \emph{Foundations
  and Trends in Communications and Information Theory}, vol.~2, no.~3, pp.
  155--239, 2006. [Online]. Available:
  \url{http://www-ee.stanford.edu/$\sim$gray/toeplitz.html}
\BIBentrySTDinterwordspacing

\bibitem{Kev01}
K.~L. Moore, ``An observation about monotonic convergence in discrete-time,
  p-type iterative learning control,'' in \emph{Proc. IEEE Int. Symposium on
  Intelligent Control}, Mexico, Sep 2001, pp. 45--49.

\bibitem{Tom89}
M.~Tomizuka, T.-C. Tsao, and K.-K. Chew, ``Analysis and synthesis of
  discrete-time repetitive controllers,'' \emph{J. Dynamic Systems, Measurement
  and Control}, vol. 111, pp. 353--358, Sep 1989.

\bibitem{Tom90}
K.-K. Chew and M.~Tomizuka, ``Digital control of repetitive errors in disk
  drive systems,'' \emph{IEEE Control Systems Mag.}, pp. 16--20, Jan 1990.

\bibitem{Won75}
B.~A. Francis and W.~M. Wonham, ``The internal model principle for linear
  multivariable regulators,'' \emph{Applied Math., Opt.,}, vol.~2, pp.
  170--194, 1975.

\bibitem{Bini}
D.~Bini and V.~Pan, ``Efficient algorithms for the evaluation of the
  eigenvalues of (block) banded toeplitz matrices,'' \emph{Mathematics of
  Computation}, vol.~50, no. 182, pp. 431--448, Apr 1988.

\bibitem{Arb}
P.~Arbenz, ``Computing eigenvalues of banded symmetric toeplitz matrices,''
  \emph{SIAM J. Sci. and Stat. Comput.}, vol.~12, no.~4, pp. 743--754, July
  1991.

\bibitem{Smi}
G.~D. Smith, \emph{Numerical Solution of Partial Differential Equations}, ser.
  Oxford Applied Mathematics and Computing Science, 2, Ed.\hskip 1em plus 0.5em
  minus 0.4em\relax Oxford, U.K.: Clarendon Press.

\end{thebibliography}
%%%%%%%%%%%%%%%%%%%%%%%%%%%%%%%%%%%%%%%%%%%%%%%%%%%%%%%%%%%%%%%%%%%%%%%%%%%%%%%%
\begin{IEEEbiography}[{\includegraphics[width=1in,height=1.25in,clip,keepaspectratio]{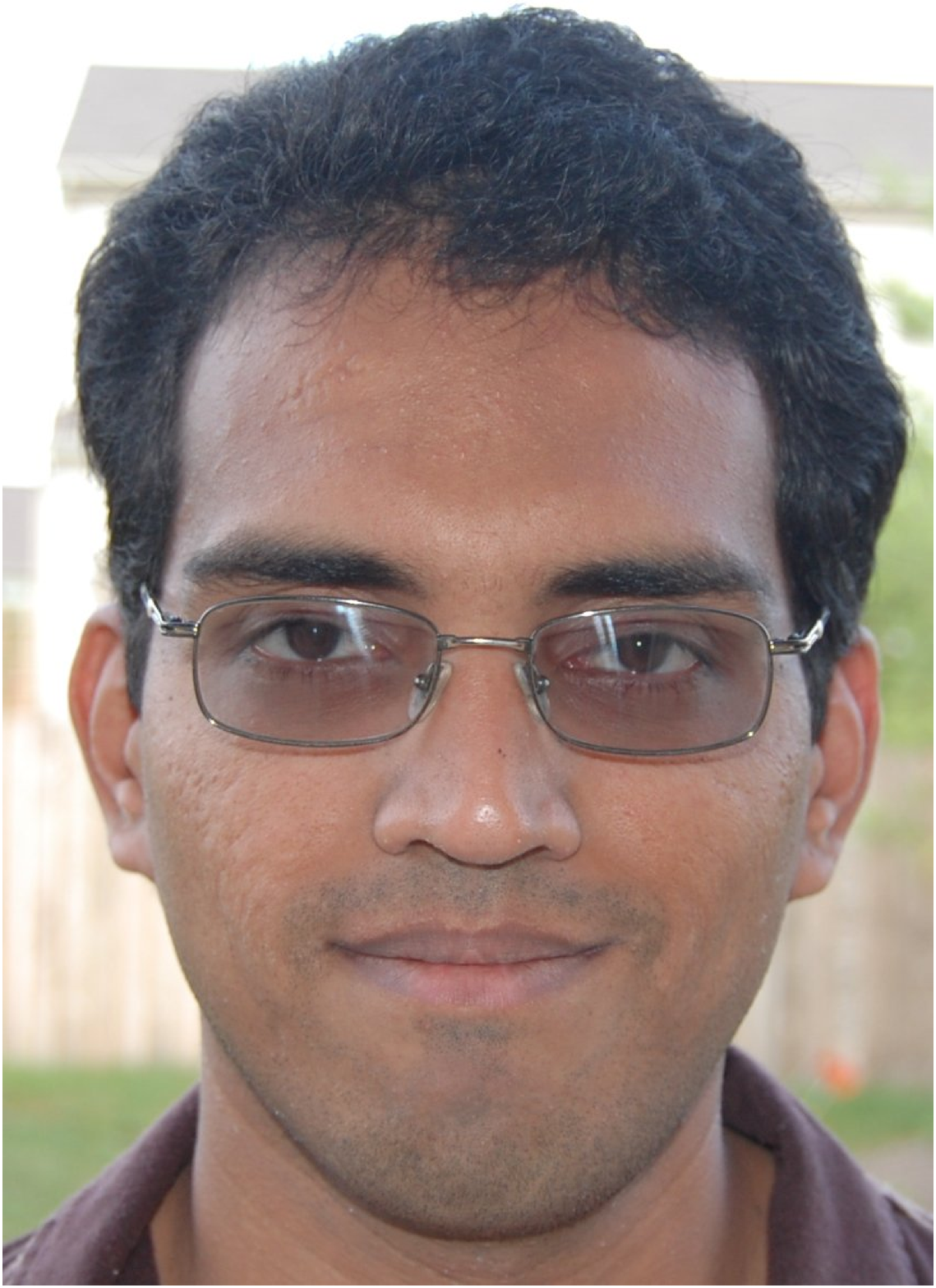}}]{Krishnamoorthy Kalyanam}
received the B.Tech. degree in mechanical engineering from
the Indian Institute of Technology Madras, Chennai, India in 2000,
and the M.S. and Ph.D. degrees in mechanical engineering from The University of California at Los Angeles, in 2003 and 2005, respectively.

In 2005 he joined the General Electric Company's Global
Research Center, Bengaluru, India as a Research Engineer.
He is currently a National Research Council Research Associate
working at the Air Force Research Lab, Wright-Patterson Air Force
Base, Ohio. His current research interests include cooperative
control of multivehicle systems, adaptive and learning control.
\end{IEEEbiography}

\begin{IEEEbiography}[{\includegraphics[width=1in,height=1.25in,clip,keepaspectratio]{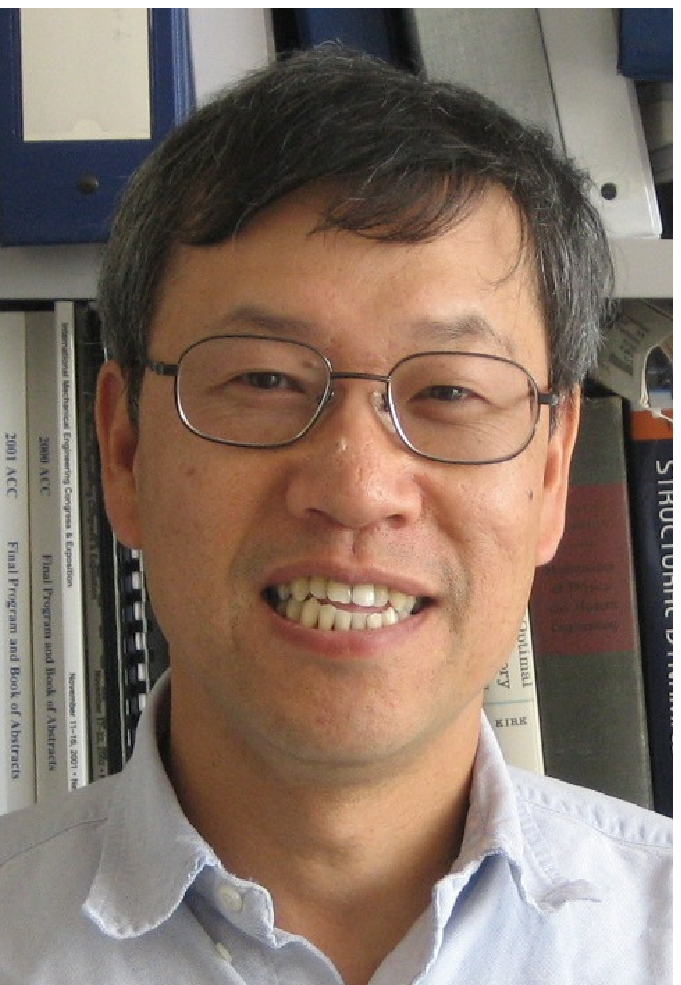}}]{Tsu-Chin Tsao}
 (S'86-M'88-SM'09) received the B.S. degree in engineering from National Taiwan University, Taipei, Taiwan, in 1981, and the M.S. and Ph.D. degrees in mechanical
engineering from The University of California at
Berkeley, in 1984 and 1988, respectively.
In 1999, he joined the faculty of the University
of California at Los Angeles, where he
is currently a Professor with the Department of Mechanical and
Aerospace Engineering. For 11
years, he was with the faculty of the Department of Mechanical and Industrial
Engineering, The University of Illinois
at Urbana-Champaign. His research interests include control systems and
mechatronics.

Prof. Tsao is a recipient of the \textit{ASME Journal of Dynamic Systems, Measurement, and Control} Best Paper Award for papers published in the journal in
1994, the Outstanding Young Investigator Award from ASME Dynamic Systems
and Control Division in 1997, and the Hugo S. Shuck Best Paper Award.
\end{IEEEbiography}

\end{document}